\documentclass[11pt,twoside]{article}
\usepackage{amsmath,amssymb,latexsym}

\newtheorem{theorem}{Theorem}[section]
\newtheorem{lemma}[theorem]{Lemma}
\newtheorem{corollary}[theorem]{Corollary}
\newtheorem{deflemma}[theorem]{Definition and Lemma}
\newcommand{\proof}{\quad{\sc Proof}:  \,}
\newcommand{\bewend}{\vspace{1ex}\hfill $ \Box $}

\newcommand{\OO}{{\cal O}}
\newcommand{\LL}{{\cal L}}
\newcommand{\BB}{{\cal B}}
\newcommand{\EE}{{\cal E}}

\newcommand{\GG}{{\cal G}}
\newcommand{\HH}{{\cal H}}
\newcommand{\A}{\mathbb A}
\renewcommand{\P}{\mathbb P}
\newcommand{\N}{\mathbb N}

\newcommand{\rank}{{\rm rk}}

\newcommand{\rel}{\mbox{$\, \omega_{X/B} \,$}}

\newcommand{\frel}{\mbox{$\, f_{\ast}\omega_{X/B} \,$}}
\newcommand{\freln}{\mbox{$\, f_{\ast}\omega_{X/B}^n \,$}}

\begin{document}

\thispagestyle{empty}

\begin{center}
{\Large \sc An Algebraic Proof of 
Iitaka's Conjecture $C_{2,1}$} 

\

by

\

{\small \sc Markus Wessler}
\end{center}

\

\begin{abstract}
We give a proof of Iitaka's conjecture $C_{2,1}$ using only elementary
methods from algebraic geometry.
\end{abstract}

\

\begin{center}
\section{\sc Introduction}
\end{center}
We consider the following situation. Let $X$ be a smooth
projective surface and let $B$ be a smooth projective curve,
defined over an algebraically closed field $k$ of characteristic
zero. Let $f : X \to B$ be a surjective morphism and let $\rel$ denote 
the relatively canonical sheaf of differentials.
Let us assume that the generic fibre is smooth of genus $g$ and let 
us denote by $\delta$ the
number of singular points in the fibres. We write
$\Lambda_n$ for the determinant of $\freln$ and $\lambda_n$
for the degree of $\Lambda_n$. Finally, let us assume that $f$ is a
relatively minimal model, which means that 
there are no exceptional curves among the
fibres.

\

In this situation, Iitaka's conjecture $C_{2,1}$ is well-known:

\

\begin{theorem}\label{iitaka}
If, in the situation above, $X_b$ denotes a general fibre, we have the 
subadditivity
$$\kappa(X) \geq \kappa(B) + \kappa(X_b)$$
of the Kodaira dimensions.
\end{theorem}

We immediately notice that, in order to give a proof for this, we 
may assume both $B$ and $X_b$ to have genus greater than zero.

\

The result \ref{iitaka} contributes to the Enriques-Kodaira
classification of surfaces, as it is presented in \cite{bpv84}, VI.,
for example.
A proof for \ref{iitaka} is given in \cite{bpv84}, III. 18.4, and as 
it is worked out there, Iitaka's conjecture
$C_{2,1}$ basically follows from: 

\begin{theorem}\label{hauptsatz1}
Keeping the assumptions made above and assuming in addition that
$f$ is non-isotrivial, we have $\lambda_1 > 0$.
\end{theorem}

There have been several kinds of proof for \ref{hauptsatz1}, 
but most of them used
analytic methods. In \cite{bpv84}, III.17, for example, a proof is 
given by considering the period map. In case this map is not 
constant, \ref{hauptsatz1} follows
from constructing a section of $\Lambda_n$, which locally arises from
modular forms. For the constant case, the image of the period map is
exactly the period point of all the smooth fibres. After excluding the
existence of singular fibres, \ref{hauptsatz1} 
then follows from Torelli's theorem.
Other proofs used analytic methods in order to achieve 
the weak positivity of 
$\frel$ (see \ref{fkpos}) like Fujita (\cite{f78}),
Kawamata (\cite{k81}) or Viehweg (\cite{v83}).  
In case algebraic methods were used, those were developed in a more
general set up in order to apply for higher
dimensions, too. 

\

But in the special case of 
a familiy of curves over a curve, one can give both an elementary and
a purely algebraic 
proof of \ref{hauptsatz1} based on positivity methods. 
The aim of this paper is to present this proof, which exclusively 
uses methods from algebraic geometry and hence yields an algebraic proof
for \ref{iitaka}.

\begin{center}
\section{\sc Positivity}
\end{center}
First we want to recall the following positivity notations and give
some properties, which we shall need below. We refer to 
\cite{ev92}, 5., for example.

\

\begin{deflemma}\label{positiv}
\rm Let $\LL$ be an invertible sheaf on a smooth projective
variety $Y$. We call $\LL$ {\it big} if its Kodaira dimension is
maximal or, equivalently, some power of $\LL$ contains an ample
subsheaf. 
We call $\LL$ {\it numerically effective} if for every
curve $C$ on $X$ one has $\deg \LL_{|_C} \geq 0$ or,
equivalently, if for every $\nu \geq 0$ and for every ample
invertible sheaf $\HH$ one has that $\LL^{\nu} \otimes \HH$ is
ample. This implies that the self intersection number of $\LL$ is
not negative and that it is positive if and only if $\LL$ is big.
Moreover, a locally free sheaf $\GG$ on $Y$ is called {\it weakly
positive} if every quotient sheaf of $\GG$ has non-negative
degree, even after a finite covering. 
\end{deflemma}

The following theorem is well known:

\begin{theorem}\label{fkpos}
The direct image sheaf $\frel$ is weakly positive.
\end{theorem}

There have been proofs by Fujita (\cite{f78}) and Kawamata
(\cite{k81}), using analytic 
methods, and
there have been several generalizations. By means of Koll\'ar's
vanishing theorem one obtains a proof (\cite{k86}, see
also \cite{v95}, Theorem 2.41), which is based on methods
from algebraic geometry.

\begin{corollary}\label{numeff}
The relatively canonical sheaf $\rel$ is numerically effective
and hence for the self intersection number we have
$c_1(\rel)^2 \geq 0$.
\end{corollary}

\proof
Let $C$ be a curve on $X$. Since $f$ is assumed to be a relatively minimal
model, we
have $c_1(\rel).C \geq 0$ in case $C$ is a fibre. If not, then the
natural map $f^{\ast}f_{\ast}\rel \to \rel$ restricted to $C$ is
surjective, and $$\deg \rel_{|_C} = c_1(\rel).C \geq 0,$$ since $\frel$
is weakly positive. 
Thus $\rel$ is numerically effective and 
by \ref{positiv} we obtain that $c_1(\rel)^2 \geq 0$.
\bewend

\

We should remark here that, in case $g > 1$, one even obtains strict 
positivity for $c_1(\rel)^2$. We will, however, not need this in the
sequel.

\begin{center}
\section{\sc Reductions}
\end{center}
In order to show \ref{hauptsatz1} we may, without loss of
generality, make some reductions, which we shall deduce in
the following. 

\

First of all, we may assume $f$ to be a semi-stable model, which means
that all the fibres are reduced normal crossing divisors. For
if this is not the case, then for some finite covering 
$\tau : B' \to B$ , for the fibre product 
$X' = X \times_B B'$ with projections $f' : X' \to B'$ and $\tau' : X'
\to X$ and for
a desingularization $d : X'' \to X'$
the induced map $f'' : X'' \to B'$ will be semi-stable.
Considering the trace map $d_{\ast}\omega_{X''} \to \omega_{X'}$,
which induces an injective map $d_{\ast}\omega_{X''/B'} 
\to \omega_{X'/B'} =
{\tau'}^{\ast}\rel$, and applying flat base change
we obtain an injective map $f''_{\ast}\omega_{X''/B'} =
f'_{\ast}d_{\ast}\omega_{X''/B'} \to f'_{\ast}{\tau'}^{\ast}\rel =
\tau^{\ast}\frel$.  Comparing degrees, we notice that if
\ref{hauptsatz1} holds for the semi-stable model $f''$, it will hold for
$f$ as well.

\

Our next reduction is based on the following
lemma which is due to Mumford (\cite{m77}, 5.10) and follows
from an easy calculation of the Riemann Roch formulae on $X$ and
on $B$.

\

\begin{lemma}\label{mumf}
For every $n \in \N$ we have
$$\lambda_n = {n \choose 2}\cdot\Big(12 \lambda_1 - \delta\Big) +
\lambda_1 =
{n \choose 2}\cdot c_1(\rel)^2 + \lambda_1$$
in case $g > 1$ and
$12 \lambda_n = n \delta = 12 n \lambda_1$
in case $g = 1$.
\end{lemma}

In case $g > 1$ the above formula tells us that 
$12 \lambda_1 \geq \delta$ (since by \ref{numeff} we have 
$c_1(\rel)^2 \geq 0$) and, moreover, $\lambda_1 = 0$
implies $\lambda_n \leq 0$ for all $n \in \N$.
In case $g = 1$ we have $12 \lambda_1 = \delta$ 
and $\lambda_n = n \cdot \lambda_1$.
Hence in both cases, in order to show $\lambda_1 > 0$, we may 
assume from now on that all the fibres are smooth and, moreover,
it suffices to show $\lambda_n > 0$ for some $n \in \N$. 

\

More reductions can be made considering the genus $g$ of the fibres.
We can immediately exclude the case $g = 1$, for if all the fibres are
smooth elliptic curves, then 
associating to each $b$ the $j$-invariant of the corresponding fibre
gives a morphism $j : B \to \A^1$, which has to be constant and
thereby forces $f$ to be isotrivial.

\

Moreover, we may assume that all the fibres are smooth
non-hyperelliptic curves of genus $g \geq 2$, as we shall prove now. To
this end, we have to distinguish two cases.
First, let us assume that  
all the fibres are smooth and
hyperelliptic curves. Then $f$ again turns out to be isotrivial. 
To see this, let us 
consider the projective bundle $\pi : \P(\frel) \to B $ associated 
to $\frel$ and let us 
denote by $\P$ the image of the morphism
$\varphi : X \to \P(\frel)$, 
which corresponds to the surjective map $f^{\ast}\frel \to
\rel$. Since all the fibres are assumed to be hyperelliptic and thus
$\varphi$ is given fibrewise by double coverings of $\P^1$, $\P$ turns
out to be a ruled suface. Denoting by $\Delta$
the discriminant, the branched covering trick (\cite{bpv84}, I.18.2)
implies that there exists an \'etale covering $\gamma : B' \to B$ such
that the pull back of $\Delta$ to $\P'$ has $2g + 2$ disjoint
components, where $\P'$ denotes the fibre product. These components
correspond to $2g + 2$ disjoint sections of the bundle $\pi' : \P' 
\to B'$, forcing it to be trivial, since $g \geq 1$. Hence the 
corresponding morphism $f' : X'
\to B'$ is isotrivial, and by property of the fibre product, so is $f$.
 
\
 
For the second case we need the following

\begin{lemma}\label{mult}
For $n \in \N$ sufficiently large, the multiplication map
$$\mu_n : S^n(\frel) \longrightarrow \freln$$
is surjective outside hyperelliptic fibres.
\end{lemma}

\proof
Let $X_b$ be a non-hyperelliptic fibre of $f$, thus $\omega_{X_b}$ is
very ample and we consider the embedding 
$X_b \to \P^{g-1}_k$ given by the global sections of $\omega_{X_b}$. We
may identify those sections with $H^0(\P^{g-1}_k ,
\OO_{\P^{g-1}_k}(1))$ which implies an isomorphism 
$$S^n(H^0(X_b , \omega_{X_b})) \to H^0(\P^{g-1}_k ,
\OO_{\P^{g-1}_k}(n))$$
for every $n \in \N$. Considering the long exact
cohomology sequence corresponding to the ideal sheaf of $X_b$ and
twisting by $n$, Serre's Vanishing Theorem implies a surjective map
$$ H^0(\P^{g-1}_k , \OO_{\P^{g-1}_k}(n)) \to 
H^0(X_b , \omega^n_{X_b}), $$
if we choose $n$ to be sufficiently large. Thus we obtain a surjective
map 
$$ S^n(H^0(X_b , \omega_{X_b})) \to H^0(X_b , \omega^n_{X_b}), $$
and by base change we are done.
\bewend

\

Now in case some of the fibres are
hyperelliptic, we consider the factorization of $\mu_n$
over its image sheaf. Since 
symmetric products of weakly positive sheaves are again weakly
positive (see for example \cite{v95}, 2.20),
$S^n(\frel)$ is weakly
positive by \ref{fkpos}, and we obtain by definition of weak
positivity that the image
sheaf has non-negative degree, hence $\lambda_n > 0$ and by \ref{mumf}
we are done. 

\

Putting all this together, it remains, in order to prove
\ref{hauptsatz1}, to show the following 

\begin{theorem}\label{hauptsatz2}
Keeping the assumptions from \ref{hauptsatz1} and assuming in addition
that $f$ is semi-stable and all the
fibres are smooth and non-hyperelliptic curves, we have $\lambda_n >
0$ for some $n \in \N$.
\end{theorem}

\begin{center}
\section{\sc The proof of \ref{hauptsatz2}}
\end{center}
In this section we prove \ref{hauptsatz2} using the
method of the universal basis as it was used by Viehweg in \cite{v89}.

\

Let $\EE$ be a locally free sheaf on $B$ of rank $m$ and let $\pi : \P
\to B$ denote the projective bundle associated to $\bigoplus^m
\EE^{\vee}$. From the construction of the projective bundle (see
for example \cite{h77}, II.7.12) we get a surjective map
$$\pi^{\ast}\bigoplus^m\EE^{\vee} \longrightarrow \OO_{\P}(1)$$
and by dualizing we obtain an injective map
$$\OO_{\P}(-1) \longrightarrow \pi^{\ast} \bigoplus^m \EE,$$
sending a local section $l$ of $\OO_{\P}(-1)$ to a tuple $(s_1 ,
\dots , s_m)$. This induces a map
$$s : \bigoplus^m \OO_{\P}(-1) \longrightarrow \pi^{\ast} \EE,$$
injective, too, sending $(f_1l , \dots , f_ml)$ to the sum of the
$f_is_i$. We call $s$ the universal basis associated to $\EE$.

\

Next we apply this method to $\EE = \frel$. Taking for $n \in \N$ 
symmetric product of the corresponding universal basis $s$
and composing with the multiplication map, we obtain a map
$$S^n\Big(\bigoplus^m \OO_{\P}(-1)\Big) \to
S^n\Big(\pi^{\ast}\frel\Big) = \pi^{\ast}S^n(\frel) \to  
\pi^{\ast}\freln,$$ which by \ref{mult} is surjective 
outside the zero divisor $D$ of $\det s$,
provided we have chosen $n$ to be sufficiently large. 
We denote the image sheaf of this morphism by $\BB$.

\

Next we consider a blowing up $\tau : \P' \to \P$ with center in $D$
such that $\BB' = \tau^{\ast}\BB$ is locally free.
Denoting by $\pi' : \P' \to B$ the composed map, we obtain a surjection
$$\rho : S^n\Big(\bigoplus^m \OO_{\P'}(-1)\Big) \longrightarrow \BB'.$$
Now since the invertible quotients
of a locally free sheaf correspond exactly to the sections of the
projective bundle (see \cite{h77}, II.7.12), we obtain
a morphism
$$\Phi' : \P' \longrightarrow \P(V),$$
where for $r =\rank \freln$ we denote by $V$ the linear space 
$\bigwedge^rS^n(k^m)$. Moreover, we have
${\Phi'}^{\ast}\OO_{\P(V)}(1) = \det \BB'$, and our next step is to
show that this sheaf is ample. This immediately follows from:

\

\begin{lemma}
In the situation above, $\Phi'$ is generically finite.
\end{lemma}

\proof
Let us denote by $G$ the Grassmannian
manifold parametrizing the $r$-dimensional quotients of the linear
space $S^n (k^m)$. Then $\Phi'$ factors 
over $G$ by considering pointwise the surjective map $\rho$, which we
recall to be given
basically by the multiplication map, and then 
composing with the Pl\"ucker embedding of $G$.
Now let us assume that there exists a generic fibre containing some 
curve $C$. Then, of
course, $C$ is not contained in any fibre of $\pi'$ over $B$, hence maps
surjectively to $B$. 
Now, considering the factorization described above, 
the image of $C$ in $G$ is assumed
to be a point, and, moreover, all the fibres $X_b$ can be
embedded in projective space by the sections of their canonical
sheaves. This implies that the global sections of all the ideal
sheaves corresponding to the fibres $X_b$, twisted by $n$, are the
same.  Since, by
increasing $n$ again if necessary, we may assume that these sheaves
are globally generated, all the fibres $X_b$ turn out to be
isomorphic, which contradicts the non-isotriviality of the family $f$.
\bewend
 
As an immediate consequence, $\Phi'$ preserves ampleness, and hence
the inclusion 

$${\Phi'}^{\ast}\OO_{\P(V)}(1)  = \det \BB' \to
{\pi'}^{\ast}\Lambda_n(nr)$$

implies that ${\pi'}^{\ast}\Lambda_n(nr)$ is big on $\P'$.

\

Let us now denote ${\pi'}^{\ast}\Lambda_n(nr)$ by $\LL$ and
the general fibre of $\pi'$ by $F$. For every $\nu \in \N$ we 
obtain the long exact cohomology sequence
$$0 \to H^0(\P',\LL^{\nu}(-F)) \to H^0(\P',\LL^{\nu}) \to 
H^0(F,\LL^{\nu}_{|_F}) \to \dots $$
Since $\LL$ is big, for $\nu$ sufficiently large $\LL^{\nu}(-F)$ 
is going to have a
non-trivial section, corresponding to a non-trivial map
$$\OO_{\P'} \longrightarrow {\pi'}^{\ast}\Lambda^{\nu}_n(\nu nr)(-F),$$
and so we obtain, for a point $P$ in general position,
$$\OO_B \longrightarrow \Lambda^{\nu}_n(-P) \otimes
\pi'_{\ast}\OO_{\P'}(\nu nr) =
\Lambda^{\nu}_n(-P) \otimes S^{\nu nr}\bigoplus^m \frel
^{\vee},$$ where the last 
equality holds by definition of the
projective bundle. Dualizing again, we get a
non-trivial map
$$S^{\nu nr}\bigoplus^m \frel \longrightarrow
\Lambda^{\nu}_n(-P).$$ Since $\frel$ is weakly positive, we
obtain by definition $\deg \Lambda^{\nu}_n(-P) \geq 0$, 
hence $\lambda_n > 0$, which completes the proof of \ref{hauptsatz2}.

\

%

\end{document}